\documentclass[10pt]{article}
\usepackage{paper,macros,graphicx}
\usepackage{epstopdf}
\usepackage{url}
\usepackage{color}
\usepackage[colorlinks=false, hidelinks]{hyperref}
\usepackage{framed}
\usepackage{subcaption}

\usepackage{algorithm}
\usepackage{algpseudocode}

\algnewcommand{\IIf}[1]{\State\algorithmicif\ #1\ \algorithmicthen}
\algnewcommand{\EndIIf}{\unskip\ \algorithmicend\ \algorithmicif}

\algnewcommand{\IfThenElse}[3]{
	\State \algorithmicif\ #1\ \algorithmicthen\ #2\ \algorithmicelse\ #3}

\usepackage[numbers, sort&compress ]{natbib}
\numberwithin{equation}{section}
\usepackage{bbm}

\usepackage{natbib}
\bibpunct[, ]{(}{)}{,}{a}{}{,}%

\usepackage{booktabs}
\usepackage{array}
\usepackage{arydshln}
\renewcommand{\thefootnote}{\arabic{footnote}}

\renewcommand{\proof}{\textbf{Proof.\ }}
\renewcommand{\qed}{~\hfill$\square$}
\renewcommand{\epsilon}{\varepsilon}

\graphicspath {{figures/}}

\title{An Adaptive Sample Size Trust-Region Method for Finite-Sum Minimization}

\author{Robert Mohr\footnotemark[1] \ \footnotemark[2]
	\and Oliver Stein\footnotemark[1]}

\begin{document}
\maketitle

\renewcommand{\thefootnote}{\fnsymbol{footnote}}
\footnotetext[1]{Institute of Operations Research, Karlsruhe Institute of Technology, Karlsruhe, Germany; E-mails: \url{robert.mohr@kit.edu}, \url{stein@kit.edu}.}
\footnotetext[2]{Corresponding author.}
\renewcommand{\thefootnote}{\arabic{footnote}}

\begin{abstract}
We propose a trust-region method for finite-sum minimization with an adaptive sample size adjustment technique, which is practical in the sense that it leads to a globally convergent method that shows strong performance empirically without the need for experimentation by the user. During the optimization process, the size of the samples is adaptively increased (or decreased) depending on the progress made on the objective function. We prove that after a finite number iterations the sample includes all points from the data set and the method becomes a full-batch trust-region method. Numerical experiments on convex and nonconvex problems support our claim that our algorithm has significant advantages compared to current state-of-the-art methods.
\end{abstract}

\section{Introduction}
\label{introduction}

There exists a variety of applications from statistics and machine learning that require the minimization of an objective function that is the sum of a large number of convex or nonconvex functions. Well known examples are logistic regression problems, the training of neural networks and nonlinear least-squares problems. In these applications, the number of functions summed up in the objective function typically corresponds to the number of data points considered. 

Successful algorithms from classical nonlinear optimization, such as quasi-Newton, nonlinear conjugate-gradient and trust-region methods, usually require the computation of the gradient and (approximate) Hessian of the objective function in every iteration. If the number of data points is very large, these computations are expensive and prohibit fast progress in the early stages of the optimization process.

Popular methods therefore use single data points or samples of data points (so-called mini-batches) in order to obtain approximate information about the objective function. Arguably the most well-known and successful algorithms in this area are the stochastic gradient descent method, which was first proposed by \cite{Robbins1951}, and its variance-reduced variants \citep[e.g.,][]{Schmidt2017, Defazio2014, Johnson2013}, in which single data points or mini-batches are used in order to approximate the gradient of the objective function. We refer the interested reader to the excellent surveys by \cite{Bottou2018} and \cite{Curtis2017b} for details concerning these methods. However, these methods have two major drawbacks. Firstly, extensive experimentation is needed for every new problem and data set in order to find hyper-parameters (e.g., the step-size) that lead to a good performance of these methods. Secondly, since only first-order information is employed, their ability to make progress in the presence of saddle points or to deal with ill-conditioned problems is limited.

In the last few years there has been growing interest in algorithms that speed up the minimization of large-scale finite-sum problems by incorporating approximate second-order information via sampling. In particular, stochastic Newton, Gauss-Newton and (limited-memory) BFGS methods were developed \citep[e.g.,][]{Byrd2011, RoostaKhorasani2019, Bollapragada2018c, Martens2010, Martens2011, Schraudolph2007, Bordes2009, SohlDickstein2014, Mokhtari2015, Byrd2016, Berahas2016, Curtis2016, Gower2016, Zhou2017, Bollapragada2018b, Berahas2019}. However, despite promising theoretical and empirical results, most of the proposed methods still depend on extensive hyper-parameter tuning for each new problem and data set. Moreover, all of these methods work with positive definite curvature approximations and experience numerical instability when these matrices become close to singular. However, in the context of nonconvex optimization, it was demonstrated by \cite{Curtis2019} that incorporating directions of negative curvature can be beneficial and, according to \cite{Dauphin2014}, they might help to escape saddle points more quickly.

In this paper, we propose a trust-region method that can be applied to large-scale nonconvex finite sum minimization. The method is very flexible with respect to the type of approximate curvature information then can be used, and can exploit directions of negative curvature. In addition to that, the method needs far less experimentation by the user, as will be demonstrated in the numerical tests.

The paper is structured as follows: After a literature review in the next section, we describe our trust-region algorithm in Section \ref{sec:ASTR}. In Section \ref{sec:theory}, we present some theoretical results concerning the convergence of the method. Afterwards, we discuss several practical considerations in Section \ref{sec:prac_Con}, before investigating the empirical performance of the method in Section \ref{sec:Num_Ex}. Section \ref{sec:final_remarks} summarizes the contributions of our paper and outlines some avenues for future research.

\section{Related Literature}

Along with nonlinear conjugate gradient and quasi-Newton methods, trust-region algorithms belong to the most reliable and efficient algorithms for the local minimization of general nonlinear functions. Theoretical results concerning global convergence properties of classical trust-region methods, as well as practical considerations, can be found in the books by \cite{Conn2000} and \cite{Nocedal2006} and the survey paper by \cite{Yuan2015}.

From a practical point of view, trust-region algorithms for the finite-sum minimization problem proposed so far can be broadly classified into three groups, depending on how sampling is used in order to obtain approximate information about the objective function. Members of the first group evaluate the objective function and its gradient exactly in each iteration, while using a sample of the data points to determine approximate curvature information \citep[e.g.,][]{Xu2019, Xu2017}. In addition to approximating curvature information, methods that belong to the second group also approximate the gradient based on a (possibly different) sample, while still evaluating the objective function exactly in every iteration \citep[e.g.,][]{Gratton2017, Yao2018, Erway2019}. The last group contains methods that, at least in the early stages of the optimization process, only work with inexact information about the objective function based on samples, i.e., the objective function is evaluated inexactly as well \citep[e.g.,][]{Chen2018, Bellavia2018, Blanchet2019}.

The main idea underlying methods from the first group is that the most expensive step in each iteration of a trust-region method is the (approximate) solution of the trust-region subproblem, at least if nontrivial curvature approximations are employed. This cost can be greatly reduced if the curvature information is approximated based on a small sample of the data points. The global convergence to first order critical points is covered by results on standard trust-region methods. However, since the objective function and its gradient are evaluated exactly in each iteration, the behavior of these methods is more similar to deterministic than to randomized methods.

This drawback also applies to the methods of the second group. In typical finite-sum problems from machine learning and statistics, the evaluation of the objective function is about half as expensive as the computation of the gradient. Therefore, although methods that approximate the gradient can be more efficient than methods that use the exact gradient, the progress of these methods will be slow in the early stages of the optimization process as long as the objective function is evaluated in every iteration.

In contrast, methods from the last group can achieve very low per iteration costs if the samples used for the approximations are sufficiently small. However, in order to obtain a convergent method, the objective function and its gradient have to be approximated with increasing accuracy. In the context of the finite-sum minimization, this necessitates increasing the corresponding sample size during the optimization process, which is referred to as dynamic/adaptive sampling or progressive batching. This technique leads to hybrid deterministic-stochastic methods, i.e., methods that start off as randomized methods and eventually turn into deterministic methods.

In these methods, the sample size can either be increased at a preset rate or adaptively according to information obtained during the optimization process. Promising theoretical and empirical results were obtained for the stochastic gradient descent and the stochastic L-BFGS methods \citep[e.g.,][]{Friedlander2012, Byrd2012, De2017, Bollapragada2018, Bollapragada2018c, Bollapragada2018b}. In the context of trust-region methods, adaptive rules for adjusting the sample size so far either depend on unknown quantities  or require experimentation by the user in order to obtain good performance.

In this paper, we describe a trust-region method for the empirical risk minimization problem with a practical adaptive sample size adjustment technique, in the sense that it leads to a globally convergent method that shows strong performance empirically without the need for experimentation by the user.

We note that the technique we propose for sample size adjustment could also be used in conjunction with the adaptive regularization method with cubics (ARC) proposed by \cite{Cartis2011} and \cite{Cartis2011b}. The ARC method is an adaptive version of the cubic regularization method first introduced by \cite{Griewank1981}. It was shown by \cite{Nesterov2006} and \cite{Cartis2011b} that cubic regularization methods and their adaptive variants are, from a worst-case complexity point of view, superior to classical trust-region methods. This fact lead to increased research interest in stochastic variants of these methods \citep[e.g.,][]{Kohler2017, Xu2017, Cartis2018}. However, we chose to propose a trust-region method since it was observed in \cite{Xu2017} that they tend to show stronger empirical performance than ARC methods.

\section{The ASTR-Algorithm}
\label{sec:ASTR}

We call our method \textbf{A}daptive \textbf{S}ample Size \textbf{T}rust-\textbf{R}egion method, or ASTR for short. It is specifically designed to solve the finite-sum minimization problem
\[
\quad \min_{x\in\mathbb{R}^d} F(x):=\frac{1}{n}\sum_{i=1}^n f_i(x),
\]
where $n,d\in\mathbb{N}$ and $f_i\in C^2(\mathbb{R}^d,\mathbb{R})$ for all $i=1,...,n$. The method consists of outer and inner iterations, shown in Algorithm \ref{alg:OuterIterations} and \ref{alg:InnerIterations}, respectively. In every inner iteration, a sample $S\subseteq\{1,...,n\}$ is chosen and Algorithm \ref{alg:tr_step} is used to compute a trust-region step for the function
$
F_S:=\frac{1}{\vert S\vert}\sum_{i\in S} f_i.
$
After a certain number of inner iterations of Algorithm \ref{alg:InnerIterations}, a candidate for the next outer iterate is returned to Algorithm \ref{alg:OuterIterations}. There, the candidate is either accepted or rejected and the sample size is adjusted. We now describe the three algorithms in detail.

\subsection{Algorithm \ref{alg:OuterIterations} - Outer Iterations}

In iteration $\nu$ of Algorithm \ref{alg:OuterIterations}, Algorithm \ref{alg:InnerIterations} is called with the current iterate $x^\nu$, sample size $s^\nu$, initial trust-region radius $\delta^\nu$ and number of inner iterations $R^\nu$ as input arguments. It returns to Algorithm \ref{alg:OuterIterations} a candidate for the next iterate $\widehat{x}^{\nu}$ and a prediction $b^{\nu}$ of the improvement in the objective function value if $\widehat{x}^{\nu}$ is accepted. Moreover, a value $\delta^{\nu+1}$ is returned, which is passed as the initial trust-region radius to Algorithm \ref{alg:InnerIterations} in the next iteration of Algorithm \ref{alg:OuterIterations}. 

The candidate $\widehat{x}^{\nu}$ is accepted if $a^{\nu}$, the  improvement in the objective function value, is nonnegative. Note that the computation of $a^{\nu}$ is an expensive operation if $n$ is large, since $F$ needs to be evaluated at $\widehat{x}^{\nu}$. The new sample size $s^{\nu+1}$ is chosen depending on the size of the ratio $\tau^\nu$ of actual to predicted improvement. A small value of $\tau^\nu$ indicates that the sampled functions used in the inner iterations do not approximate $F$ accurately enough and that the sample size should therefore be increased. A large value of $\tau^\nu$, however, is an indicator that faster progress in the inner iterations might be possible if the sample size is decreased. 

Note that every time the sample size is increased/decreased in the outer iteration, the inner iterations get more computationally expensive/cheap and the number of inner iterations should therefore also be decreased/increased. In Section \ref{sec:prac_Con} we explain how to update the sample size and the number of inner iterations in order to obtain a method with strong empirical performance.

\begin{algorithm}[tb]
	\caption{ASTR - Outer Iterations}
	\label{alg:OuterIterations}
	\begin{algorithmic}[1]
		\State {\bfseries Input:} Initial point $x^0\in\mathbb{R}^d$, parameters $s^0,\widehat{R}\in\mathbb{N}$ and $\delta^0, \theta>0$;
		\For{$\nu=0,1,...$} 
		\State Select the number of inner iterations $R^\nu\in\{1,...,\widehat{R}\}$;
		\State Compute $\widehat{x}^{\nu}$, $b^{\nu}$ and $\delta^{\nu+1}$ via Algorithm \ref{alg:InnerIterations} with inputs $x^{\nu}, s^\nu$, $\delta^{\nu}, R^\nu$;
		\If{$s^\nu<n$}
		\State Set $a^{\nu} = F(x^{\nu} )- F(\widehat{x}^{\nu})$;
		\IfThenElse{$a^{\nu} \geq 0$}{set $x^{\nu+1} = \widehat{x}^{\nu}$}{set $x^{\nu+1} = x^{\nu}$};
		\IfThenElse{$b^{\nu} > 0$}{set $\tau^\nu = a^{\nu}/b^{\nu}$}{set $\tau^\nu= 0$};
		\IfThenElse{$\tau^\nu < \theta$}{select $s^{\nu+1} \in\{s^\nu+1,...,n\}$}{ $s^{\nu+1}\in\{1,...,s^\nu\}$};
		\Else
		\State Set $x^{\nu+1} = \widehat{x}^{\nu}$ and $s^{\nu+1}=n$;
		\EndIf
		\EndFor
	\end{algorithmic}
\end{algorithm}

\subsection{Algorithm \ref{alg:InnerIterations} - Inner Iterations}

In iteration $k$ of Algorithm \ref{alg:InnerIterations}, a sample $S^{\nu,k}$ of size $s^\nu$ is chosen. For notational convenience, we define 
$
F^{\nu,k}:=F_{S^{\nu,k}}
$
and
$
g^{\nu,k}:=\nabla F^{\nu,k}(x^{\nu,k}).
$
If the gradient $g^{\nu,k}$ does not satisfy $\Vert g^{\nu,k}\Vert_2\geq\epsilon$,  where $\epsilon$ is a preset threshold, no step is taken and a new sample is selected in the next iteration. 
Otherwise, a trust-region step $d^{\nu,k}$ is computed and the initial trust-region radius $\delta^{\nu,k+1}$ for the next iteration is determined with Algorithm \ref{alg:tr_step}.

After $R$ iterations, the last inner iterate $x^{\nu,R}$ and the current trust-region radius $\delta^{\nu,R}$ are returned to Algorithm \ref{alg:OuterIterations}. Additionally, the average improvement on the sampled functions during the inner iterations
$$
b^{\nu} := \frac{1}{R} \sum_{k=0}^{R-1} b^{\nu,k} = \frac{1}{R} \sum_{k=0}^{R-1} \big(F^{\nu,k}(x^{\nu,k})- F^{\nu,k}(x^{\nu,k+1})\big)
$$
is returned to Algorithm \ref{alg:OuterIterations} as a prediction for the improvement on the objective function $F$ if $x^{\nu,R}$ is accepted as the next outer iterate.

\begin{algorithm}[tb]
	\caption{ASTR - Inner Iterations}
	\label{alg:InnerIterations}
	\begin{algorithmic}[1]
		\State {\bfseries Input:} $x^\nu, \delta^\nu, s^\nu$ and $R^\nu$ from Algorithm \ref{alg:OuterIterations} and parameter $\epsilon >0$;
		\State {Set $x^{\nu,0} = x^\nu$, $\delta^{\nu,0}=\delta^\nu$ and $R=R^{\nu}$;}
		\For{$k=0,1,...,R-1$} 
		\State Choose a sample $S^{\nu,k}\subseteq\{1,...,n\}$ of size $s^\nu$;
		\If{$\Vert g^{\nu,k}\Vert_2\geq\epsilon$}
		\State {Compute $d^{\nu,k}, \delta^{\nu,k+1}$ via Algorithm \ref{alg:tr_step} with inputs $x^{\nu,k}$, $g^{\nu,k}, \delta^{\nu,k}$};
		\State Set $x^{\nu,k+1} = x^{\nu,k} +d^{\nu,k}$;
		\State Set $b^{\nu,k}=F^{\nu,k}(x^{\nu,k})- F^{\nu,k}(x^{\nu,k+1})$;
		
		\Else
		\State Set $x^{\nu,k+1} = x^{\nu,k}$, $b^{\nu,k}=0$ and $\delta^{\nu,k+1}= \delta^{\nu,k}$;
		\EndIf
		\EndFor
		\State {\bfseries Output:} $\widehat{x}^{\nu} = x^{\nu,R}$, $b^{\nu} = \displaystyle\frac{1}{R}\sum_{k=0}^{R-1} b^{\nu,k}$, $\delta^{\nu+1}=\delta^{\nu,R}$;
	\end{algorithmic}
\end{algorithm}

\subsection{Algorithm \ref{alg:tr_step} - Trust-region step}

In iteration $r$ of Algorithm 3 an (approximate) solution $d^r$ to the trust-region subproblem
\begin{align}
\label{prob:trs}
\min_{d\in\mathbb{R}^d}\ m^{\nu,k}(d) \quad \text{s.t.} \quad \Vert d\Vert_2\leq \Delta_r,
\end{align}
is computed, where $m^{\nu,k}$ is a (quadratic) model of $F^{\nu,k}$ at $x^{\nu,k}$ defined as
\[
m^{\nu,k}(d) := F^{\nu,k}(x^{\nu,k})+\langle g^{\nu,k},d \rangle + \frac{1}{2}d^\intercal A^{\nu,k} d,
\]
and $\Delta_r$ is the current trust-region radius. The matrix $A^{\nu,k}$ can be used to include curvature information in the model. However, it is also possible to only use first-order information by setting $A^{\nu,k}=0$. 

An (approximate) solution $d^{r}$ to problem \eqref{prob:trs} is accepted if the ratio
\begin{align}
\label{tr_ratio}
\rho^{\nu,k}(d^r):=\displaystyle\frac{F^{\nu,k}(x^{\nu,k} +d^r)-F^{\nu,k}(x^{\nu,k})}{m^{\nu,k}(0)-m^{\nu,k}(d^{r})}
\end{align}
of the actual improvement on $F^{\nu,k}$ to the improvement predicted by the quadratic model is above a certain threshold $\eta_1$. The intuition behind this is that if the ratio $\rho^{\nu,k}(d^r)$ is above this threshold, this is an indication that the model $m^{\nu,k}$ is a good approximation of $F^{\nu,k}$ on the feasible set of problem \eqref{prob:trs} (the so called ``trust-region''). As long as the ratio \eqref{tr_ratio} is smaller then $\eta_1$, the trust-region radius is decreased and a new (approximate) solution of the trust-region subproblem is computed. Since we use the exact gradient of $F^{\nu,k}$ in our model $m^{\nu,k}$, it is always possible to find an acceptable (approximate) solution to the trust-region subproblem if the trust-region radius is sufficiently small (see Theorem \ref{the:accepted_if_small_tr}). Note that if the ratio \eqref{tr_ratio} is not only larger than $\eta_1$ but also larger than $\eta_2$, then the trust-region radius is increased such that larger steps might be taken in the next inner iteration.

\begin{algorithm}[tb]
	\caption{Trust-region step}
	\label{alg:tr_step}
	\begin{algorithmic}[1]
		\State {\bfseries Input:} $x^{\nu,k}$, $g^{\nu,k}$ and $\delta^{\nu,k}$ from Algorithm \ref{alg:InnerIterations} and parameters $\eta_1, \eta_2\in(0,1)$ and $0<\gamma_1<1<\gamma_2$;
		\State Set $r= -1$, $\Delta_0 = \delta^{\nu,k}$;
		\Repeat 
		\State $r \leftarrow r+1$;
		\State Compute an (approximate) solution $d^{r}$ of problem
		\begin{align*}
		\min_{d\in\mathbb{R}^d}\ m^{\nu,k}(d) \quad \text{s.t.} \quad \Vert d\Vert_2\leq \Delta_r;
		\end{align*}
		\State Compute $\rho^{\nu,k}(d^r):=\displaystyle\frac{F^{\nu,k}(x^{\nu,k} +d^r)-F^{\nu,k}(x^{\nu,k})}{m^{\nu,k}(0)-m^{\nu,k}(d^{r})}$;
		\State Set $\Delta_{r+1}=\gamma_1 \Vert d^r\Vert_2$;
		\Until{$\rho^{\nu,k}(d^r)<\eta_1$}
		
		\IfThenElse{$\rho^{\nu,k}(d^r)\geq\eta_2$ and $\Vert d^r\Vert_2 = \Delta_r$}{set $\Delta_+ = \gamma_2\Delta_r$}{set $\Delta_+ = \Delta_r$};
		\State {\bfseries Output:} $d^{\nu,k}=d^r$ and  $\delta^{\nu,k+1}=\Delta_+$;
	\end{algorithmic}
\end{algorithm}

\section{Theoretical Analysis}
\label{sec:theory}

In this section we prove that after a finite number of iterations, the sample size $s^\nu$ reaches $n$ and the ASTR method becomes a full-batch trust-region method.

\subsection{Assumptions}

\begin{assumption}
	\label{ass_low_bounded}
	The function $F$ is bounded below on $\mathbb{R}^d$, i.e., there exists a constant $\kappa_1\in\mathbb{R}$ such that
	$
	F(x)\geq \kappa_1
	$
	for all $x\in\mathbb{R}^d$.
\end{assumption}

\begin{assumption}
	\label{ass_lipschitz}
	The functions $f_i$, $i=1,...,n$, are twice continuously differentiable and their gradients $\nabla f_i$ are Lipschitz continuous.
\end{assumption}

Note that Assumption \ref{ass_lipschitz} implies that the Hessians $D^2f_i(x)$ are uniformly bounded in $x$ for all $i$. From the triangular inequality it immediately follows that the Hessian of the function $F_S$ is uniformly bounded in $x$ and $S$, i.e., there exists a constant $\kappa_2>0$ such that the inequality
\begin{align}
\label{hes_bounded}
\Vert D^2F_S(x)\Vert_2 \leq \kappa_2
\end{align}
holds for any $S\subseteq\{1,...,n\}$ and $x\in\mathbb{R}^d$.

\begin{assumption}
	\label{ass_Curvature}
	There exists a constant $\kappa_{3}>0$ such that for all $\nu, k$ we have that 
	$
	\Vert A^{\nu,k}\Vert_2 \leq \kappa_{3}.
	$
\end{assumption}

Assumption \ref{ass_Curvature} is trivially satisfied if $A^{\nu,k}=0$ for all $\nu, k$. Due to \eqref{hes_bounded} we know that Assumption \ref{ass_Curvature} is also satisfied if we set $A^{\nu,k} = D^2F_S(x)$ for any $S$ and $x$.

\begin{assumption}
	\label{ass_Cauchy}
	There exists a constant $\kappa_{4}\in(0,1)$ such that for all $\nu, k, r$ we have that 
	\[
	m^{\nu,k}(0) - m^{\nu,k}(d^r)\geq \kappa_{4} \Vert g^{\nu,k}\Vert_2\min(\frac{\Vert g^{\nu,k} \Vert_2}{1+\Vert A^{\nu,k}\Vert_2},\Delta_r).
	\]
\end{assumption}

If $d^r$ is a sufficiently accurate approximation of the exact solution of the trust-region subproblem, Assumption \ref{ass_Cauchy} is satisfied. Moreover, there exists a variety of methods for the inexact solution of the trust-region subproblem such that Assumption \ref{ass_Cauchy} is satisfied, e.g., the truncated conjugate gradient (CG) method by \cite{Toint1981} and \cite{Steihaug1983} or the truncated Lanczos method by \cite{Gould1999}.

\subsection{Theoretical Results}

The first two theorems presented in this section are modifications of well known results in the literature on trust-region methods, see e.g. Theorems 6.4.2 and 6.4.3 in \cite{Conn2000}. The proofs of these two theorems are provided in the appendix.

\begin{theorem}
	\label{the:accepted_if_small_tr}
	Suppose A\ref{ass_lipschitz}, A\ref{ass_Curvature} and A\ref{ass_Cauchy} hold. Then there exists a constant $\kappa_5>0$ such that for all $\nu$, $k$ and $r$ the inequality $\Vert d^r\Vert_2\leq\kappa_5$ implies that $\rho^{\nu,k}(d^r)\geq\eta_1$ holds.
\end{theorem}

The previous theorem guarantees that Algorithm \ref{alg:tr_step} terminates after a finite number of steps. Moreover, it is instrumental in proving the following result.

\begin{theorem}
	\label{the:tr_not_too_small}
	Suppose A\ref{ass_lipschitz}, A\ref{ass_Curvature} and A\ref{ass_Cauchy} hold. Then there exists a constant $\kappa_{6}>0$ such that $\delta^{\nu,k}\geq \kappa_{6}$ holds for all $\nu,k$.
\end{theorem}

\begin{theorem}
	\label{the:b_bounded}
	Suppose A\ref{ass_lipschitz}, A\ref{ass_Curvature} and A\ref{ass_Cauchy} hold. Then there exists a constant $\kappa_{7}>0$ such that for all $\nu$ either $b^\nu = 0$ or $b^\nu \geq\kappa_{7}$ holds.
\end{theorem}

\proof For any $k$, if $\Vert g^{\nu,k}\Vert_2< \epsilon$, then $b^{\nu,k} = 0$. On the other hand, if $\Vert g^{\nu,k}\Vert_2\geq \epsilon$, then Theorem \ref{the:accepted_if_small_tr} guarantees that Algorithm \ref{alg:tr_step} terminates with a trust-region step $d^{\nu,k}$ that satisfies $\rho^{\nu,k}(d^{\nu,k})\geq\eta_1$. Thus, we obtain from \eqref{tr_ratio} and Assumption \ref{ass_Cauchy} that
\begin{align*}
b^{\nu,k} &= F^{\nu,k}(x^{\nu,k})- F^{\nu,k}(x^{\nu,k+1})
\geq \eta_1 \left(m^{\nu,k}(0)-m^{\nu,k}(d^{\nu,k})\right)\\
&\geq \eta_1\kappa_{4} \Vert g^{\nu,k}\Vert_2\min(\frac{\Vert  g^{\nu,k} \Vert_2}{1+\Vert A^{\nu,k}\Vert_2},\Delta^{\nu,k}),
\end{align*}
where $\Delta^{\nu,k}$ denotes the trust-region radius in the  subproblem that was used to compute the trust-region step $d^{\nu,k}$. From Theorem \ref{the:tr_not_too_small} it follows that $$\Delta^{\nu,k} \geq \frac{\delta^{\nu,k}}{\gamma_2} \geq \frac{\kappa_6}{\gamma_2} $$
and from Assumption \ref{ass_Curvature} we know that $$\frac{\Vert  g^{\nu,k} \Vert_2}{1+\Vert A^{\nu,k}\Vert_2} \geq \frac{\Vert  g^{\nu,k} \Vert_2}{1+\kappa_3}.$$ We therefore have that
\[
b^{\nu,k}\geq\eta_1\kappa_{4} \epsilon\min(\frac{\epsilon}{1+\kappa_{3}},\frac{\kappa_6}{\gamma_2}).
\]
Thus, since $b^\nu =\displaystyle \frac{1}{R^\nu}\sum_{k=0}^{R^\nu-1} b^{\nu,k}$, we either have $b^\nu = 0$ or $b^\nu \geq\kappa_{7}$, with
\[
\kappa_7:=\frac{\eta_1\kappa_{4} \epsilon}{\widehat{R}}\min(\frac{\epsilon}{1+\kappa_{3}},\frac{\kappa_6}{\gamma_2})>0.
\]\qed

Concerning the previous theorem, we note that $b^\nu = 0$ can only occur if $\Vert g^{\nu,k}\Vert_2<\epsilon$ for all inner iterations.

\begin{theorem}
	\label{the:subset_size}
	Suppose A\ref{ass_low_bounded}, A\ref{ass_lipschitz}, A\ref{ass_Curvature} and A\ref{ass_Cauchy} hold. Then there exists a $\nu_0\in\mathbb{N}$ such that $s^\nu = n$ holds for all $\nu\geq \nu_0$.
\end{theorem}

\proof We show that there exists a $\nu_0\in\mathbb{N}$ such that $s^{\nu_0} = n$ since this implies the assertion in the theorem. Assume, for the purpose of deriving a contradiction, that $s^\nu < n$ for all $\nu\in\mathbb{N}$. Since $\tau^\nu<\theta$ implies $s^{\nu+1}>s^\nu$, there does not exist a $\bar{\nu}\in\mathbb{N}$ such that $\tau^\nu<\theta$ for all $\nu>\bar{\nu}$. Thus, there exists a subsequence $(x^{\nu(j)})$ with $\tau^{\nu(j)}\geq\theta$ for all $j\in\mathbb{N}$. This implies that $b^{\nu(j)} >0$ for all $j\in\mathbb{N}$. From Theorem \ref{the:b_bounded} we now obtain that $b^{\nu(j)} >\kappa_{7}$ and therefore
\begin{align*}
F(x^{\nu(j)} )- F(\widehat{x}^{\nu(j)})=a^{\nu(j)}\geq\theta b^{\nu(j)}
\geq \theta \kappa_{7}
\end{align*} 
for all $j\in\mathbb{N}$. Since the sequence $(F(x^\nu))$ is monotonically nonincreasing we obtain for all $i\in\mathbb{N}$
\begin{align*}
F(x^0)-F(x^{\nu(i)+1})
&= \sum_{\nu=0}^{\nu(i)} (F(x^\nu)-F(x^{\nu+1}))
\geq \sum_{j=0}^i (F(x^{\nu(j)})-F(x^{\nu(j)+1}))\\
&= \sum_{j=0}^i (F(x^{\nu(j)})-F(\widehat{x}^{\nu(j)}))
\geq (i+1)\theta \kappa_{7},
\end{align*}
and therefore
\begin{align*}
F(x^0)-F(x^{\nu(i)+1})\overset{i\rightarrow\infty}{\rightarrow} +\infty,
\end{align*}
in contradiction to Assumption \ref{ass_low_bounded}.\qed

Theorem \ref{the:subset_size} ensures that after a finite number of iterations, the ASTR method becomes a standard (full-batch) trust-region method. It implies that global convergence of the ASTR method follows from the global convergence results about standard trust-region methods, e.g., Theorem 6.4.6 in \cite{Conn2000}.

\section{Practical Considerations}
\label{sec:prac_Con}

In the description of our algorithm in Section \ref{sec:ASTR} we left out several details that do not need specification in order to prove the theoretical results in Section \ref{sec:theory}, but which are nonetheless important with regard to the practical implementation of the method. The purpose of this section is to close this gap.

\subsection{Random Sampling and Sample Size Adjustment}
\label{sample_size_adj}

We propose to select the samples in the inner iterations of our algorithm uniformly at random, although other selection strategies (deterministic and stochastic) are possible and may be worth investigating. 

\cite{Friedlander2012} and \cite{Byrd2012} showed that when the sample size is increased geometrically in stochastic gradient decent, then the expected optimality gap converges linearly . Inspired by this strategy, we propose to choose a constant $\omega>1$ and set
$
s^{\nu+1} = \min(\lceil\omega s^{\nu} \rceil,n)
$
whenever $\tau^\nu < \theta$. If $\tau^\nu \geq \theta$, one can simply set $s^{\nu+1} = s^\nu$, i.e., the sample size is never decreased. We leave the question of whether strategies for decreasing the sample size can lead to performance benefits for future research.

\subsection{Incorporation of Curvature Information and the Solution of the Trust-Region Subproblems}

For $A^{\nu,k}=0$ the solution of the trust-region subproblem (\ref{prob:trs}) is
$
d^r = -\Delta_r g^{\nu,k}/\Vert g^{\nu,k} \Vert_2.
$
Thus, if $A^{\nu,k}=0$ for all $\nu,k$, the ASTR method is a adaptive sample size gradient method. 

It is one of the strengths of the ASTR method that any kind of curvature information can be used. However, one has to keep in mind that the choice of $A^{\nu,k}$ is tightly coupled with effort necessary to compute an (approximate) solution to the trust-region subproblem. Fortunately, the trust-region subproblem is a problem that has been studied for decades and one can choose from a wide variety of methods in order to determine exact or inexact solutions, see, for example, \cite{Conn2000}. Consequently, there is a lot of flexibility for investigating different ways to incorporate curvature information.

The most straightforward way to incorporate curvature information is to set $A^{\nu,k}=D^2 F^{\nu,k}(x^{\nu,k})$ as soon as the sample size $s^{\nu}$ is considered large enough for the sampled Hessian to contain meaningful curvature information. If the decision variable is very high dimensional, the (sampled) Hessian of the objective is expensive to compute and might be too large to store. However, if the truncated CG method is used for the solution of the trust-region subproblems, only matrix-vector products of $A^{\nu,k}$ and certain vectors need to be computed. These matrix-vector products can be efficiently computed for various problems in supervised machine learning without ever forming the matrix $A^{\nu,k}$ explicitly, see \cite{Pearlmutter1994}. This technique is known as ``Hessian-free'' optimization. 

Note that if this ``Hessian-free'' technique is used, the cost of multiplying $A^{\nu,k}$ with a vector depends on the size of the sample used for the computation of $A^{\nu,k}$. This cost can therefore be reduced by setting $A^{\nu,k}=D^2 F_{S_H^{\nu,k}}(x^{\nu,k})$ for some subsample $S_H^{\nu,k}\subseteq S^{\nu,k}$, provided that $s^{\nu}$ is large enough. In \cite{Xu2017} one can find some guidance on how to subsample the Hessian in a trust-region framework when exact gradient information is used. 

Also note that for nonconvex problems, $A^{\nu,k}$ can be indefinite and directions of negative curvature can be exploited. This might be particularly useful in the proximity of saddle points, which are considered one of the main obstacles when training neural networks with current methods, see, for example, \cite{Dauphin2014}. 

\subsection{Adjusting the Number of Inner Iterations}
\label{sec:Numb_Inner_Iter}

The last detail that needs to be specified is how the number of inner iterations $R^{\nu+1}$ should be adjusted depending on the updated sample size $s^{\nu+1}$. We suggest to select $R^{\nu+1}$ in a way such that the total computational cost in all the inner iterations combined is approximately equal to the cost of evaluating the objective function $F$ in the outer iteration. Consequently, the number of inner iterations depends on the kind of curvature information used and the method for the solution of the trust-region subproblem. 

To be more concrete: Assume that evaluating the objective function is half as expensive as computing the gradient, and that the computation of a Hessian-vector product costs approximately the same as the computation of a gradient. This is indeed the case for various applications, see Section \ref{sec:Num_Ex} for some examples. 

If $A^{\nu,k}=0$ for some $\nu$ and all $k$, the solution of the trust-region subproblem \eqref{prob:trs} is given explicitly by $d^r = -\Delta_r g^{\nu,k}/\Vert g^{\nu,k} \Vert_2$. Thus, the cost of one inner iteration corresponds to the cost of evaluating the gradient $g^{\nu, k}$. Consequently, $R^{\nu}$ should satisfy the equation
$
R^{\nu} \cdot s^{\nu}/n = 0.5
$
and we obtain $R^{\nu} = n/(2s^{\nu})$ for the number of inner iterations of Algorithm \ref{alg:InnerIterations}.

If $A^{\nu,k}=D^2 F_{S_H^{\nu,k}}(x^{\nu,k})$ for some $\nu$ and all $k$, and the truncated CG method is used to approximately solve the trust-region subproblems, analogous reasoning leads to the formula
$
R^{\nu} = n/((2 + \bar{\alpha}) \cdot s^{\nu} + \bar{\beta}\cdot 2 \cdot  s_H^\nu),
$
where $\bar{\alpha}$ denotes the average number of iterations of Algorithm 3, $\bar{\beta}$ denotes the average number of iterations the truncated CG method requires to find an approximate solution to the trust-region subproblems and $s_H^\nu:=\vert S_H^{\nu,k}\vert$, where we assume that the size of the subsample $S_H^{\nu,k}$ is fixed during the inner iterations.

\section{Numerical Experiments}
\label{sec:Num_Ex}

In this section we compare the ASTR method with a mini-batch stochastic gradient descent method (SGD), the SVRG method by \cite{Johnson2013} and a full-batch Trust-Region Newton-CG method (TR). We consider three classification problems: logistic regression (convex), nonlinear least-squares (nonconvex) and neural network training (nonconvex). For each problem and data set, we report the training errors of the methods against CPU time measurements. The algorithms were implemented in Python and the computations were performed on an Intel Core i7-9700K with 32 GB of main memory.

In contrast to the SGD and SVRG methods, who depend on hyper-parameter tuning for reasonable performance, we did not perform hyper-parameter tuning to individual problems or data sets for the ASTR or TR methods, i.e., we always used the same hyper-parameters. 

\textbf{TR:} For the standard Trust-Region Newton-CG method, we used $\delta^0=1$ as the initial trust-region radius and the standard parameters from the literature for the trail point acceptance and trust-region radius update, see \cite{Conn2000}.  The maximum number of conjugate-gradient iterations was set to 30. 

\textbf{ASTR:} For the parameters that the ASTR and TR methods have in common, we used the same values. For the additional parameters, we chose $\theta = 0.5$ and $s^0 = 0.01\cdot n$. We increased the sample size as described in Section \ref{sample_size_adj} with $\omega=2$. The parameter $\epsilon$ was set close to machine precision. For the curvature information in the ASTR method, we chose  $A^{\nu,k}:=D^2 F_{S_H^{\nu,k}}(x^{\nu,k})$, where $S_H^{\nu,k}\subseteq S^{\nu,k}$ with $s_H^{\nu}:=\vert S_H^{\nu,k} \vert = 0.1\cdot s^{\nu}$ and adjusted the number of inner iterations as it was described in section \ref{sec:Numb_Inner_Iter} (with $\bar{\alpha} = 5$ and $\bar{\beta} = 20$). When $s^{\nu}$ reaches its maximal size of $n$, $s_H^{\nu}$ is doubled in every outer iteration until it reaches $n$ as well.

\textbf{SGD:} Two hyper-parameters where tuned for each problem and data set. The best combination of a step-size $t\in\{10^{-6},10^{-5} \ldots, 1, 10\}$ and mini-batch size $s=\lceil\zeta \cdot n\rceil$ for $\zeta\in\{\frac{1}{n}, 10^{-5}, 10^{-4} \ldots, 10^{-1}\}$ was selected.

\textbf{SVRG:} For each problem and data set we tried all combinations of step-sizes $t\in\{10^{-6},10^{-5} \ldots, 1, 10\}$ and number of inner iterations $K=\lceil\mu \cdot n\rceil$ for $\mu\in\{10^{-4}, 10^{-3}, \ldots, 1, 2\}$, and selected the combination that achieved the minimal training error. 

For each problem and data set, a starting point $x^0$ was randomly generated and used by each of the methods. We ran each method for a fixed time budget. This was also the time budget that SGD and SVRG were tuned to. In order to report the training error $F(x^k) - F^\star$, the value $F^\star$ was determined by running the full-batch TR-Algorithm until it was unable to improve the objective value due to numeric precision.

\begin{table}[H]
	\caption{Data sets used in the numerical experiments.If no split of the data set into training and test set was provided in the source, we chose 10\% of the data points randomly as our test set. Otherwise we kept the original split into training and test points, except for the data set \emph{ijcnn}, where we shrank the size of the test set from 65\% to 10\% of the data points.}
	\begin{tabular}{cccccc}
		\hline
		Data set & \# training/test points  & \# Features & \# Classes & Source \\ \hline
		a9a & 32,561/16,281 & 123 & 2 & \cite{Platt1998} &  \\
		w8a & 49,749/14,951 & 300 & 2 & \cite{Platt1998} &  \\
		odd\_even & 60,000/10,000 & 784 & 2 & \cite{Lecun1998b} &  \\
		ijcnn & 127,522/14,169 & 22 & 2 & \cite{Chang2011} & \\
		skin & 220,551/24,506 & 3 & 2 & \cite{Chang2011} &  \\
		covertype & 522,911/58,101 & 54 & 2 & \cite{Collobert2002} &  \\
		SUSY & 4,500,000/500,000 & 18 & 2 & \cite{Baldi2014} &  \\
		HIGGS & 9,900,000/1,100,000 & 28 & 2 & \cite{Baldi2014} &  \\
		MNIST & 60,000/10,000 & 784 &10 & \cite{Lecun1998b} & \\ \hline
	\end{tabular}
	\label{tab:datasets}
\end{table}

\subsection{Logistic Regression}

Given a data set $(z^i,y_i)\in\mathbb{R}^d\times\{-1,1\}$, $i=1,...,n$, we consider the $\ell_2$-regularized logistic regression problem
\begin{align*}
\min_x\  F(x)=\frac{1}{n}\sum_{i=1}^n \log\big(1+ e^{-y_i(x^\intercal z_i)}\big) + \lambda \Vert x\Vert_2^2\quad \text{with}\ \lambda = \frac{1}{n}.
\end{align*}
Since the objective $F$ is strongly convex, there exists a globally minimal point and it coincides with the unique critical point of $F$. We test the methods on the binary classification data sets described in Table \ref{tab:datasets}. 

\begin{figure}
	\centerline{\includegraphics[width=\textwidth]{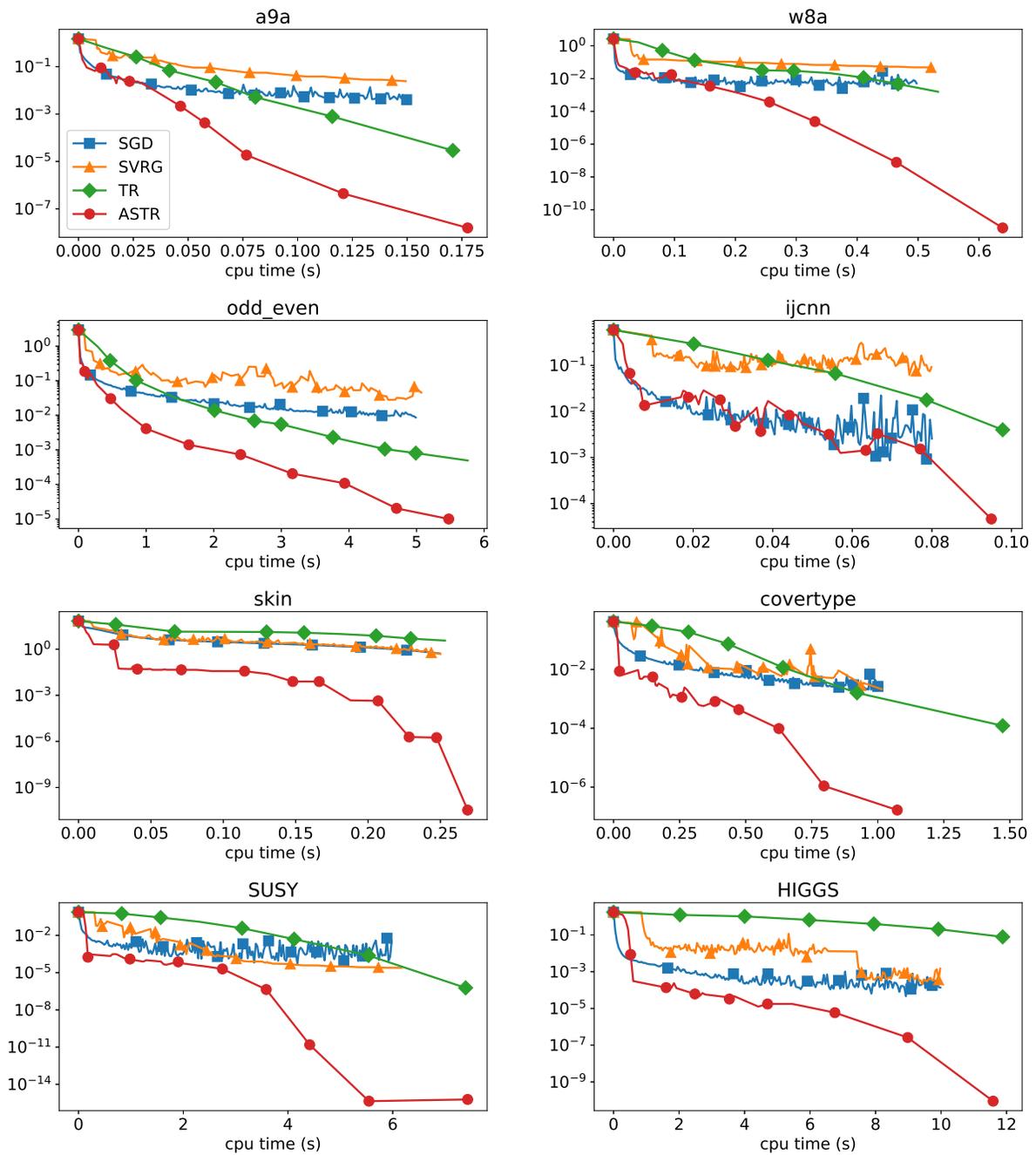}}
	\caption{Logistic regression error of SGD, SVRG, TR and ASTR on the different data sets.\label{fig:LR_Training_Errors}}
\end{figure}

In Figure \ref{fig:LR_Training_Errors} we report the results concerning the minimization of the training error. We observe consistent superior performance of the ASTR method compared to the full-batch TR and the tuned SGD and SVRG methods. 

\begin{figure}
	\centerline{\includegraphics[width=\textwidth]{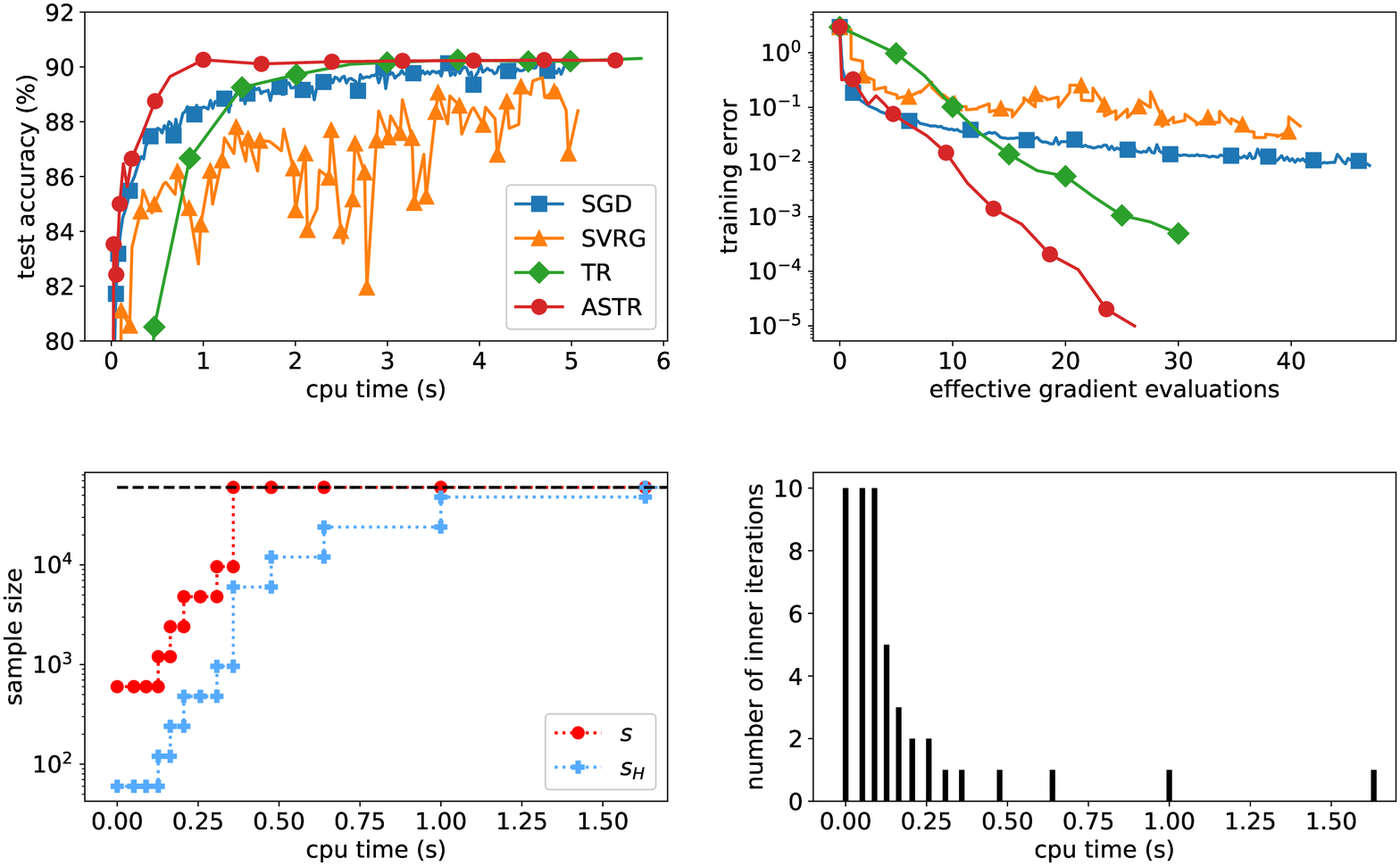}}
	\caption{Further details on the performance of SGD, SVRG, TR and ASTR on the \emph{odd\_even} data set. In the plots in the second row, the behaviour of the sample sizes and the number of inner iterations of ASTR is depicted.\label{fig:LR_odd_even}}
\end{figure}

\begin{figure}
	\centerline{\includegraphics[width=\textwidth]{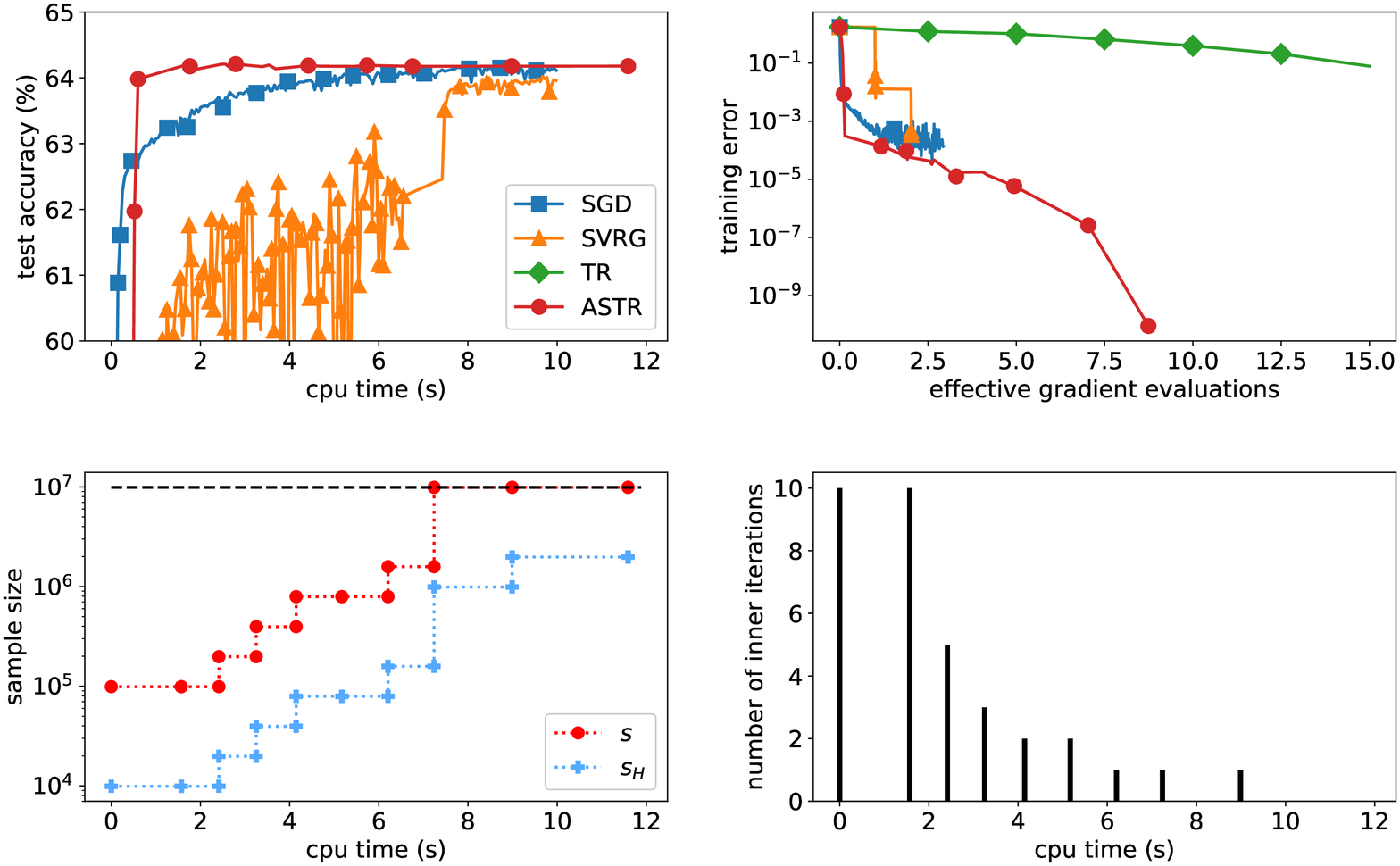}}
	\caption{Further details on the performance of SGD, SVRG, TR and ASTR on the \emph{HIGGS} data set. In the plots in the second row, the behaviour of the sample sizes and the number of inner iterations of ASTR is depicted.\label{fig:LR_HIGGS}}
\end{figure}

We note that ASTR and the tuned SGD method consistently need less CPU time than the other methods in order to archive a high test accuracy. ASTR is faster than SGD for the data sets \emph{odd\_even}, \emph{skin}, \emph{covertype} and \emph{HIGGS}, and equally fast for all the remaining data sets. The test accuracies for \emph{odd\_even} and \emph{HIGGS} are shown in the upper left panels of Figures \ref{fig:LR_odd_even} and \ref{fig:LR_HIGGS}, respectively.

We also compared the performance of the algorithms with respect to \emph{effective gradient evaluations}, a platform and implementation independent measure often used in the literature, see, for example, \cite{Bollapragada2018b}. In order to determine the number of effective gradient evaluations per iteration for each method, we made use of the fact that for each of the problems considered in this section, evaluating the objective is half as expensive as computing its gradient. And the latter operation costs the same as computing a Hessian-vector product. As to be expected, the qualitative results of the comparison of the algorithms remains unchanged when this alternative measure is used, see, for example, the upper right panels of Figures \ref{fig:LR_odd_even} and \ref{fig:LR_HIGGS}, where the training errors for the data sets \emph{odd\_even} and \emph{HIGGS} are depicted. However, since the methods we compare are very dissimilar, we believe that CPU times are more transparent and therefore more appropriate in order to evaluate the performance of the methods.

Finally, the plots in the last row of Figures \ref{fig:LR_odd_even} and \ref{fig:LR_HIGGS} provide additional details concerning the ADST method. On the left, we report the behavior of the sample sizes used in the gradient and Hessian matrix approximations ($s$ and $s_H$, respectively), on the right, the corresponding number of inner iterations is depicted. 

\subsection{Nonlinear Least-Squares}

We now focus on binary classification with squared loss as a concrete instance of a nonlinear (and nonconvex) least-squares problem. Given a data set $(z^i,y_i)\in\mathbb{R}^d\times\{0,1\}$, $i=1,...,n$, we consider the problem
\begin{align*}
\min_x\  F(x)=\frac{1}{n}\sum_{i=1}^n \big( y_i - \phi(x^\intercal z^i)\big)^2,
\end{align*}
where $\phi$ denotes the sigmoid function, i.e., $\phi(t) =  \frac{1}{1+e^{-t}}$. We use the same data sets that were used for logistic regression in the previous section.

\begin{figure}
	\centerline{\includegraphics[width=\textwidth]{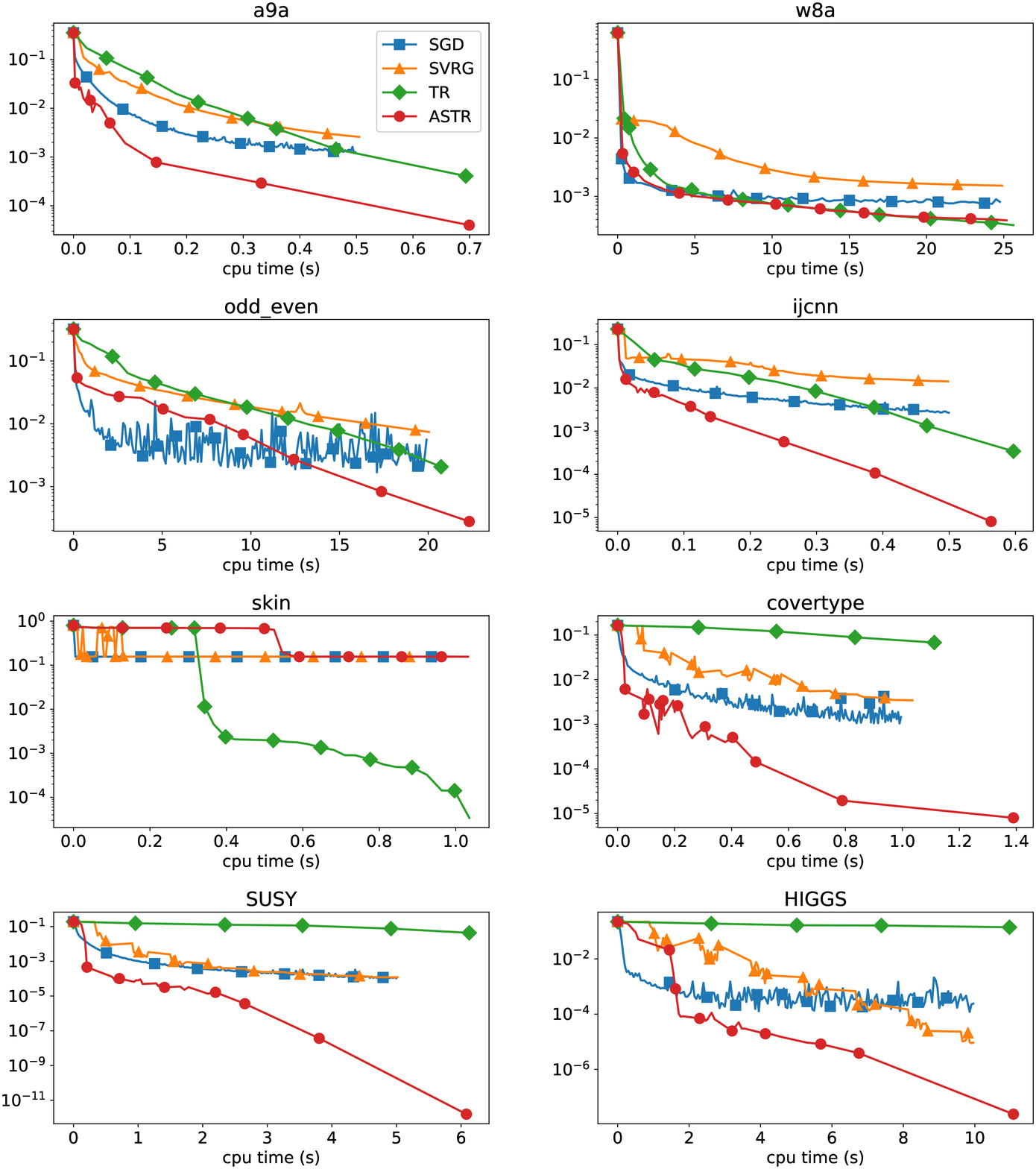}}
	\caption{Nonlinear least-squares error of SGD, SVRG, TR and ASTR on the different data sets.\label{fig:NLS_Training_Errors}}
\end{figure}

In Figure \ref{fig:NLS_Training_Errors} we report the results of the numerical experiments. Again, our results show strong performance of the ASTR method. On all data sets, except for \emph{w8a} and \emph{skin}, ASTR is clearly superior to the other methods concerning the minimization of the training error. For the data set \emph{skin}, it seems like ASTR, SGD and SVRG approximate a local minimal point, whereas TR approximates either a better local or the actual global minimal point.

Again, we point out that the test accuracies of ASTR and the tuned SGD method are comparable and that they are superior to the SVRG and TR methods. Only on the data set \emph{skin}, the TR method approximates a minimal point with much better generalization properties than the local minimal point approximated by the other methods.

\subsection{Neural Network Training}

Finally, we also considered the problem of training a simple two layer feed-forward neural network on the popular \emph{MNIST} data set of handwritten digits, see \cite{Lecun1998b}. The fully connected two layer neural network has 748 input neurons, 100 hidden neurons and 10 output neurons. The hidden neurons implement the logistic function, the output neurons the softmax function. Thus, if we have a data set $(z^i,y_i)\in\mathbb{R}^{784}\times\{0,1\}^{10}$, $i=1,...,n$, and choose the cross-entropy loss function, we arrive at the optimization problem
\[
\min_{x\in\mathbb{R}^d} \  F(x):=-\sum_{i=1}^n \sum_{j=1}^{10} y^i_j\ln( [h(z^i,x)]_j),
\]
where $h(\cdot, x)$ denotes the function that implements the neural network with weight vector $x$. 

\begin{figure}
	\centerline{\includegraphics[width=\textwidth]{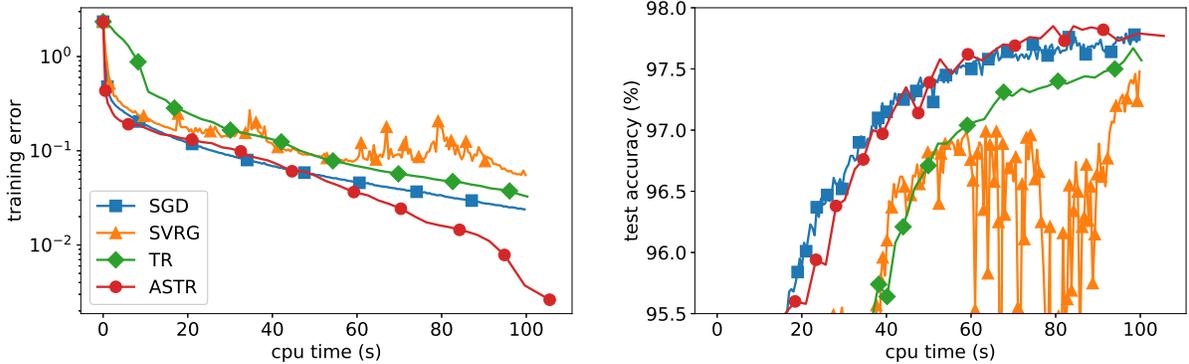}}
	\caption{Training error and test accuracy of SGD, SVRG, TR and ASTR on the MNIST dataset.\label{fig:NN_Training_Error_Accuracy}}
\end{figure}

In Figure \ref{fig:NN_Training_Error_Accuracy} one can observe that ASTR archives better results than the other methods concerning the training error. With regard to the test accuracy, ASTR performes on par with the tuned SGD method.

\section{Final Remarks}
\label{sec:final_remarks}

In this paper, we proposed an adaptively sampled trust-region method for finite-sum minimization. We showed theoretically that the sample size is eventually increased to the size of the whole training data set, which implies global convergence. Our numerical experiments demonstrated strong performance of our method, which did not involve any hyper-parameter tuning to individual problems or data sets.

A promising avenue for future research is the incorporation of different kinds of curvature information. Limited-memory techniques could be used to make our method applicable to high-dimensional problems, see \cite{Burdakov2017} and \cite{Erway2019}. Approximations to the diagonal of the Hessian, as described in \cite{Gower2018}, could make the method more efficient for very large data sets and the training of more complex neural networks.

\bibliographystyle{abbrvnat}
\bibliography{ASTR_Paper_ArXiv}

\appendix
\clearpage
\section{Proofs of Theorem \ref{the:accepted_if_small_tr} and Theorem \ref{the:tr_not_too_small}}
\label{appendix_proofs}

\subsection{Proof of Theorem \ref{the:accepted_if_small_tr}}

Define the constant
\[
\kappa_5:=\frac{(1-\eta_1)\kappa_{4}\epsilon}{1+\kappa_2+\kappa_{3}}
\]
and assume that $\Vert d^r\Vert_2\leq\kappa_5$ holds. From the definition of $\rho^{\nu,k}(d^r)$ it follows that
\begin{align}
1-\rho^{\nu,k}(d^r) 
&=  1-\frac{F^{\nu,k}(x^{\nu,k})- F^{\nu,k}(x^{\nu,k} +d )}{m^{\nu,k}(0)-m^{\nu,k}(d^r)}\nonumber\\
\label{eq:ratio}
&= \frac{F^{\nu,k}(x^{\nu,k} +d )-m^{\nu,k}(d^r)}{m^{\nu,k}(0)-m^{\nu,k}(d^r)},
\end{align}
where we used that $m^{\nu,k}(0)=F^{\nu,k}(x^{\nu,k})$. From Assumption \ref{ass_Cauchy} we know that
\[
m^{\nu,k}(0) - m^{\nu,k}(d^r)\geq \kappa_{4} \Vert g^{\nu,k}\Vert_2\min(\frac{\Vert g^{\nu,k} \Vert_2}{1+\Vert A^{\nu,k}\Vert_2},\Delta_r),
\]
where $\Delta_r$ denotes the trust-region radius in the  subproblem that was used to compute the trust-region step $d^r$.
With $\Delta_r\geq\Vert d^r\Vert_2$, $\Vert g^{\nu,k}\Vert_2\geq \epsilon$ and Assumption \ref{ass_Curvature} we obtain
\[
m^{\nu,k}(0) - m^{\nu,k}(d^r) \geq \kappa_{4}\epsilon\min(\frac{\epsilon}{1+\kappa_{3}},\Vert d^r\Vert_2).
\]
Moreover, due to $\eta_1,\kappa_{4}\in(0,1)$ and $\kappa_2\geq 0$ we know that $\Vert d^r\Vert_2\leq\kappa_5$ implies $\Vert d^r\Vert_2\leq\frac{\epsilon}{1+\kappa_{3}}$. Thus, we have
\[
m^{\nu,k}(0) - m^{\nu,k}(d^r)\geq \kappa_{4} \epsilon \Vert d^r\Vert_2.
\]
This inequality together with \eqref{eq:ratio} yields
\begin{align*}
1-\rho^{\nu,k}(d^r)
&\leq \frac{F^{\nu,k}(x^{\nu,k} +d^r )-m^{\nu,k}(d^r)}{\kappa_{4} \epsilon \Vert d^r\Vert_2 }.
\end{align*}

Now, it follows from Taylor's theorem that for some $\lambda$ in the line segment $[x^{\nu,k},x^{\nu,k} +d]$ it holds that
\[
F^{\nu,k}(x^{\nu,k}+d^r) = F^{\nu,k}(x^{\nu,k}) + \langle g^{\nu,k},d^r \rangle + \frac{1}{2}(d^r)^\intercal D^2 F^{\nu,k}(\lambda) d^r
\]
and together with Assumptions \ref{ass_lipschitz} and \ref{ass_Curvature} we obtain
\begin{align*}
F^{\nu,k}(x^{\nu,k}+d^r)-m^{\nu,k}(d^r)
&= \frac{1}{2}(d^r)^\intercal D^2 F^{\nu,k}(\lambda) d^r -  \frac{1}{2}(d^r)^\intercal A^{\nu,k} d^r\\
&\leq \frac{1}{2}\left( \Vert D^2 F^{\nu,k}(\lambda) \Vert_2 + \Vert A^{\nu,k} \Vert_2 \right)\Vert d^r\Vert_2^2\\
&\leq (\kappa_2 + \kappa_{3})\Vert d^r\Vert_2^2.
\end{align*}
Thus, we have that
\begin{align*}
1-\rho^{\nu,k}(d^r) 
&\leq \frac{(\kappa_2 + \kappa_{3})}{\kappa_{4} \epsilon }\Vert d^r\Vert_2
\leq \frac{(\kappa_2 + \kappa_{3})}{\kappa_{4} \epsilon }\kappa_5
= \frac{(\kappa_2 + \kappa_{3})}{\kappa_{4} \epsilon }\frac{(1-\eta_1)\kappa_{4}\epsilon}{1+\kappa_2+\kappa_{3}}\\
&= \frac{(\kappa_2 + \kappa_{3})}{1+\kappa_2+\kappa_{3}}(1-\eta_1)
\leq 1-\eta_1,
\end{align*}
and therefore $\rho^{\nu,k}(d^r)\geq\eta_1$.\qed

\subsection{Proof of Theorem \ref{the:tr_not_too_small}}

Define the constant
\[
\kappa_{6}:=\min(\frac{\gamma_1(1-\eta_1)\kappa_{4}\epsilon}{1+\kappa_2+\kappa_{3}},\delta^{0,0}).
\]

Assume, for the purpose of deriving a contradiction, that $(\nu,k)$ is the first iteration such that $\delta^{\nu,k}< \kappa_{6}$.
Since
$
\kappa_{6}\leq \delta^{0,0}
$
we have $(\nu,k)\neq (0,0)$. Moreover, due to $\delta^{\nu-1,R}=\delta^{\nu,0}$ we have $k\geq 1$. The value $\delta^{\nu,k}$ is calculated via Algorithm $\ref{alg:tr_step}$ with $\Delta_0=\delta^{\nu,k-1}$ as the initial trust-region. Since $(\nu,k)$ is the first iteration such that $\delta^{\nu,k}< \kappa_{6}$, we know that  $\Delta_0=\delta^{\nu,k-1}\geq \kappa_{6}$. Since $\delta^{\nu,k}=\Delta_+$ and $\Delta_+\geq \Delta_r$ it follows that $\Delta_r<\kappa_{6}$. Thus, there must exist a smallest index $j\in\{1,...,r\}$ such that $\Delta_j<\kappa_{6}$. Clearly, if $j$ is the first iteration such that $\Delta_j<\kappa_{6}$ holds, it must hold that $\rho^{\nu,k}(d^{j-1})<\eta_1$.
Consequently, we have that
\[
\Delta_j=\gamma_1\Vert d^{j-1}\Vert_2
\]
and thus
\begin{align*}
\Vert d^{j-1}\Vert_2
= \frac{\Delta_j}{\gamma_1}
\leq \frac{\kappa_{6}}{\gamma_1}
\leq \frac{(1-\eta_1)\kappa_{4}\epsilon}{1+\kappa_2+\kappa_{3}}
= \kappa_5,
\end{align*}
with $\kappa_5$ as defined in the proof of Theorem \ref{the:accepted_if_small_tr}.
From Theorem \ref{the:accepted_if_small_tr} it now follows that the inequality $\rho^{\nu,k}(d^{j-1})\geq\eta_1$ holds, which is a contradiction since we already argued that $\rho^{\nu,k}(d^{j-1})<\eta_1$.\qed

\end{document}